\documentclass[a4paper,12pt]{article}

\usepackage[utf8]{inputenc}                  
\usepackage[T1]{fontenc}                     
\usepackage{amsthm,amsmath,amssymb}          
\usepackage{mathtools}                       
\usepackage{mathrsfs}                        
\usepackage{dsfont}                          
\usepackage{hyperref}                        

\newcommand{\ud}{\mathrm{d}}
\newcommand{\e}{\mathrm{e}}
\newcommand{\CR}{\mathds{R}}
\newcommand{\CZ}{\mathds{Z}}
\newcommand{\CC}{\mathds{C}}

\newcommand{\dep}[2]{\displaystyle\frac{\partial\emph{$#1$}}{\partial{\emph{$#2$}}}}
\newcommand{\depp}[2]{\textstyle\frac{\partial\emph{$#1$}}{\partial{\emph{$#2$}}}}

\newcommand{\escalar}[2]{\left\langle{\emph{$#1$}},{\emph{$#2$}}\right\rangle}

\newcommand{\rrightarrow}{\mathrel{\mathrlap{\rightarrow}\mkern1mu\rightarrow}}

\theoremstyle{plain} 
\newtheorem{theorem}{Theorem}[section]

\newtheorem{conjecture}{Conjecture}[section]
\theoremstyle{definition}

\newtheorem{remark}{Remark}[section]

\begin{document}


\title{Symplectic structures and globally hyperbolic spacetimes}
\author{Romero Solha\footnote{
Address: Avil\'{e}s, Spain (33401). Email: romerosolha@gmail.com}}
\date{\today}

\maketitle


\begin{abstract}
The aim of this note is to present a construction of symplectic structures on orientable globally hyperbolic $4$-dimensional lorentzian manifolds. Said structures are defined on the manifold itself, not on its cotangent bundle. It also includes a discussion about their geometric quantisation.
\end{abstract}



\section{Introduction}

\hspace{1.5em}It is well known that the cotangent bundle of any manifold admits a symplectic structure. It is also a straightforward corollary of the results in \cite{Giroux} to obtain a symplectic structure on any orientable globally hyperbolic 4-dimensional lorentzian manifold (theorem \ref{contactstructure}). However, these constructions do not use (or are related to) a lorentzian structure.

Here a symplectic structure $\varpi$ on a globally hyperbolic manifold $(M,\mathrm{g})$ is obtained via an orthogonal decomposition of its tangent bundle into two plane subbundles, $TM=N_\star\oplus N$. The bundle $N$ is trivial and it is the normal bundle to a symplectic foliation of $M$, whilst $N_\star$ is the tangent bundle to said foliation.

For each of those subbundles of $TM$ a complex structure is constructed (which in turn makes $TM$ itself a complex bundle), together with hermitian inner products and connexions. The sum of the real parts of these hermitian inner products equals the euclidean structure obtained via a Wick rotation of $\mathrm{g}$, and the connexions are related to projections of the Levi-Civita connexion to each subbundle, with $\varpi$ being the sum of the imaginary parts of the curvatures.

For the Schwarzschild solution $\varpi$ is shown to be a symplectic structure, in general it is conjectured that this happens when $N_\star$ is not trivial (conjecture \ref{symplecticstructure}) for a specific choice of connexions on $N$. The closedness is guaranteed by the fact that it is a sum of curvatures of complex line bundles, the nondegeneracy seems to be controlled by $N_\star$, for Minkowski and de Sitter spacetimes both do not support nontrivial complex line bundles and $\varpi=0$.

Since a symplectic structure is provided for the spacetime itself, one can apply the machinery of geometric quantisation to $(M,\varpi)$. Due to $\varpi$ being already the sum curvatures of $N_\star$ and $N$, the associated prequantum line bundle will be $N_\star\otimes N$, and because $M$ can be understood as the configuration space (not as a phase space), there is no reason to eliminate degrees of freedom via polarisations; thus, the Hilbert space taken is the one formed by square integrable (with respect to $\varpi\wedge\varpi$) sections of the bundle $N_\star\otimes N$. In subsection \ref{quantumgravity} it is suggested how to interpret the quantisation of $(M,\varpi)$. Subsection \ref{symplecticfoliationsection} elaborates on the aforementioned symplectic foliation and its associated Poisson structure, which could also be exploited for a deformation quantisation approach.

Schwarzschild's solution is used in section \ref{schwarzschild} as an example to illustrate the constructions.

Unless otherwise stated, all objects will be smooth; manifolds are real, Hausdorff, paracompact, and connected (finite dimensionality of some homology and cohomology groups might be implicitly assumed). Physical quantities are measured using the International System of Units.


\subsection{Acknowledgements}

\hspace{1.5em}A relationship between newtonian gravity and the existence of an underlying symplectic foliation was recognised by the author during a research stay at UPV/EHU, whilst working under the Basque government project IT1094-16; where the author learned about cosymplectic structures by collaborating with Maria Fern\'{a}ndez, Eva Miranda, and C\'{e}dric Oms. The extension of the Reeb vector field to a $4$-dimensional globally hyperbolic manifold was recognised by the author through a series of meetings with Miguel S\'{a}nchez during a research stay at UGR. The author also wants to thank for a comment made by Jonathan Weitsman on a first attempt at the construction of a symplectic structure for the Schwarzschild solution. Weitsman observed that what was done was closer to quantising the cotangent bundle with the vertical polarisation, since the end result also gives square integrable functions on $M$, but the fact that the prequantum line bundle was not flat on the base suggested that it was similar to a twisted cotangent bundle.


\section{Orientable globally hyperbolic manifold}

\hspace{1.5em}If $(M,\mathrm{g})$ is a lorentzian manifold, then a $1$-dimensional submanifold $\gamma\subset M$ is a causal curve when the lorentzian structure $\mathrm{g}$ restricted to its tangent bundle (wherever it is defined\footnote{Often, a causal curve is only continuous with an open and dense subset of it being smooth.}) is nonpositive; in case it is strictly negative it is called a timelike curve, and in case it vanishes everywhere it is called a lightlike curve or null curve.

A parametrised curve from an open subset of $\CR$ to $M$ is unextendible whenever it does not attain any limits over the bounday of its domain.

On a lorentzian manifold $(M,\mathrm{g})$, a subset $Q\subset M$ is a Cauchy hypersurface if all unextendible causal curves in $M$ intersects $Q$, and no more than once when they are timelike. A lorentzian manifold $(M,\mathrm{g})$ admitting a Cauchy hypersurface is called a globally hyperbolic manifold \cite{Wald}.

\begin{theorem}[cf. \cite{BeralSanchez}]\label{hyperbolicspacetime}
For a globally hyperbolic spacetime $(M,\mathrm{g})$, there exists a pair of functions $\phi,t\in C^\infty(M)$ such that $t:M\to\CR$ defines a fibration which the fibres are isometric to a 3-dimensional Cauchy hypersurface $Q:=t^{-1}(\{0\})$, and $(M,\mathrm{g})$ is isometric to $(Q\times\CR,\mathrm{g}_t-c^2\e^{2\phi}\ud t\otimes\ud t)$, where for each $a\in\CR$ one has $\mathrm{g}_a:=\mathrm{g}|_{t^{-1}(\{a\})}$.
\end{theorem}

\begin{remark}
Assuming $M$ orientable, if $M\cong Q\times\CR$, then $Q$ is orientable, and one can apply Whitehead's theorem \cite{Whitehead} to conclude that the tangent bundle of $M$ is trivial.
\end{remark}

Before even attempting to define a symplectic structure over globally hyperbolic spacetimes, it is worth noticing that such structures actually exist.

\begin{theorem}\label{contactstructure}
Any orientable $4$-dimensional globally hyperbolic manifold $(M,\mathrm{g})$ can be endowed with a symplectic structure.
\end{theorem}
\begin{proof}
Since there exists a hypersurface $Q\subset M$ such that $M$ is diffeomorphic to $Q\times\CR$, the orientability and dimension of $M$ implies that $Q$ is an orientable $3$-dimensional manifold; wherefore, \cite{Giroux} guarantees that $Q$ admitts a contact structure, and taking a symplectisation of it provides a symplectic structure for $Q\times\CR$.
\end{proof}

It might be possible to achieve a meaningful contact structure that captures the essencial gravitational aspects of an orientable $4$-dimensional globally hyperbolic spacetime. Unfortunately, instead of constructing a contact structure on $Q$ from the $1$-parameter family of riemannian structures induced by $\mathrm{g}$, another aproach is followed in this text. Indeed, the relevant $1$-form is further from a contact form, since it defines an integrable distribution, as oppose to a maximally nonintegrable one.


\subsection{Future and gravitational vector fields}

\hspace{1.5em}Under the hypothesis of theorem \ref{hyperbolicspacetime}, $M$ is endowed with a nonvanishing exact $1$-form $\ud t\in\Omega^1(M)$. The kernel of $\ud t$ defines a distribution $\mathrm{ker}(\ud t)\subset\mathfrak{X}(M)$, as well as an associated hyperplane subbundle of $TM$, for which its integrability is guaranteed by the exactness of $\ud t$, and its leaves are exactly the (Cauchy hypersurfaces) fibres of theorem \ref{hyperbolicspacetime}.

It is also possible to choose a vector field $\depp{}{t}\in\mathfrak{X}(M)$ satisfying $\ud t(\depp{}{t})=1$ (e.g. by taking an arbitrary riemannian structure on $M$ and normalising the dual to $\ud t$); hence, the $C^\infty(M)$-module of vector fields (orthogonally) decomposes as
\begin{equation*}
\mathfrak{X}(M)=\mathrm{ker}(\ud t)\oplus\left\langle\dep{}{t}\right\rangle \ .
\end{equation*}

Considering the Levi-Civita connexion associated to $\mathrm{g}$,
\begin{equation*}
\nabla:\mathfrak{X}(M)\to\Omega^1(M)\otimes\mathfrak{X}(M) \ ,  
\end{equation*}since $\nabla\mathrm{g}=0$, one can define
\begin{equation*}
X:=-\frac{1}{c\e^\phi}\dep{}{t}\in\mathfrak{X}(M)
\end{equation*}and
\begin{equation*}
\tau:=\imath_X\mathrm{g}\in\Omega^1(M)
\end{equation*}such that
\begin{equation*}
\tau=c\e^\phi\ud t \ ,
\end{equation*}
\begin{equation*}
\tau(X)=\mathrm{g}(X,X)=-1 \ ,
\end{equation*}and
\begin{equation*}
0=X(-1)=X(\mathrm{g}(X,X))=\mathrm{g}(\nabla_XX,X)+\mathrm{g}(X,\nabla_XX)=2\mathrm{g}(\nabla_XX,X) \ ;
\end{equation*}thus,
\begin{equation*}
R:=-c^2\nabla_XX\in\mathrm{ker}(\ud t) \ ,
\end{equation*}and
\begin{equation*}
\theta:=\imath_R\mathrm{g}\in\Omega^1(M)
\end{equation*}satisfies
\begin{equation*}
\theta(X)=-c^2\mathrm{g}(\nabla_XX,X)=0 \ ,
\end{equation*}i.e.
\begin{equation*}
X\in\mathrm{ker}(\theta) \ .
\end{equation*}

From the definition of the Hodge star operator for $1$-forms,
\begin{equation*}
\theta\wedge\star\tau=\mathrm{g}(R,X)\star 1=0 \ .
\end{equation*}It also holds true that
\begin{equation*}
\theta=-c^2\imath_X\ud\tau \ .
\end{equation*}Using $\nabla\mathrm{g}=0$ and $\mathrm{tor}(\nabla)=0$, for any pair $Z_1,Z_2\in\mathfrak{X}(M)$,
\begin{equation*}
\ud\tau(Z_1,Z_2)=\mathrm{g}(\nabla_{Z_1}X,Z_2)-\mathrm{g}(\nabla_{Z_2}X,Z_1) \ ;
\end{equation*}hence,
\begin{equation*}
\imath_X\ud\tau(Z_2)=\mathrm{g}(\nabla_XX,Z_2)-\mathrm{g}(\nabla_{Z_2}X,X) \ ,
\end{equation*}and because
\begin{equation*}
0=Z_2(-1)=Z_2(\mathrm{g}(X,X))=\mathrm{g}(\nabla_{Z_2}X,X)+\mathrm{g}(X,\nabla_{Z_2}X)=2\mathrm{g}(\nabla_{Z_2}X,X) \ ,
\end{equation*}one has
\begin{equation*}
\imath_X\ud\tau(Z_2)=-\frac{1}{c^2}\mathrm{g}(R,Z_2)=-\frac{1}{c^2}\theta(Z_2) \ .
\end{equation*}Also,
\begin{equation*}
\ud\tau=\ud\phi\wedge\tau
\end{equation*}implies
\begin{equation*}
\theta=c^2\dep{\phi}{t}\ud t-c^2\ud\phi \ ;
\end{equation*}therefore, $R$ is the gradient vector field of $-c^2\phi$ with respect to the riemannian structure on $Q$ (the restriction of $\mathrm{g}$ to $\mathrm{ker}(\ud t)$),
\begin{equation*}
\tau\wedge\theta=c^2\ud\tau \ ,
\end{equation*}
\begin{equation*}
\ud\theta=-c^2\ud t\wedge\ud\left(\dep{\phi}{t}\right) \ ,
\end{equation*}and
\begin{equation*}
\ud t\wedge\ud\theta=0
\end{equation*}guarantees that $\ud\theta$ vanishes when restricted to $\mathrm{ker}(\ud t)$.


\section{Complex structure}

\hspace{1.5em}An oriented $4$-dimensional globally hyperbolic manifold (possibly excluding some of its points) can be endowed with a complex structure, and its tangent bundle can be decomposed as the sum of two complex line bundles.

\begin{remark}
What some authors call an almost complex strucure is called here a complex structure. The reason being that it makes the tangent bundle a complex bundle, whilst when it is integrable it provides a holomorphic bundle, and it seems more consistent to call such a structure a holomorphic one (as oppose to other authors who call it an integrable almost complex structure).
\end{remark}

The complement of $\mathrm{ker}(\theta)\cap\mathrm{ker}(\ud t)$ is generated by $R$ and $X$, since $\mathrm{ker}(\ud t)$ restricts the vector fields to the ones tangent to $Q$ and $\mathrm{ker}(\theta)$ is the orthogonal complement of $R$; hence,
\begin{equation*}
\mathfrak{X}(M)=\mathrm{ker}(\theta)\cap\mathrm{ker}(\ud t)\oplus\langle R\rangle\oplus\langle X\rangle \ .
\end{equation*}This orthogonal decomposition of the $C^\infty(M)$-module of vector fields comes with an orthogonal projection
\begin{equation*}
\mathds{1}-\frac{\theta\otimes R}{\mathrm{g}(R,R)}+\tau\otimes X:\mathfrak{X}(M)\rrightarrow\mathrm{ker}(\theta)\cap\mathrm{ker}(\ud t) \ .
\end{equation*}


\subsection{Foliation complex line bundle}

\hspace{1.5em}Because $\ud\theta$ vanishes when restricted to $\mathrm{ker}(\ud t)$, the intersection of the kernels of $\theta$ and $\ud t$ defines an integrable distribution; where $R$ does not vanish, its leaves are $2$-dimensional. Ergo, (still where $R$ does not vanish) the subbundle of $TM$ associated to $\mathrm{ker}(\theta)\cap\mathrm{ker}(\ud t)$ is a plane subbundle, and a complex structure on each of its fibres is defined by the unique solution $\mathrm{j}_\star \in\mathrm{Aut}(\mathrm{ker}(\theta)\cap\mathrm{ker}(\ud t))$ of
\begin{equation*}
\mathrm{g}\circ(\mathrm{j}_\star \oplus\mathds{1})=-\frac{\imath_R\star\tau}{\sqrt{\mathrm{g}(R,R)}} \ .
\end{equation*}It is convenient to extend it to the whole of $\mathrm{End}(\mathfrak{X}(M))$ as
\begin{equation*}
\mathrm{j}_\star R=\mathrm{j}_\star X=0 \ .
\end{equation*}

Denoting by $N_\star$ the complex line bundle created from $\mathrm{ker}(\theta)\cap\mathrm{ker}(\ud t)$ and $\mathrm{j}_\star $, a hermitian inner product is defined by
\begin{equation*}
\mathrm{h}_\star:=\mathrm{g}+\sqrt{-1}\mathrm{g}\circ(\mathrm{j}_\star \oplus\mathds{1}) \ ,
\end{equation*}and the projection of the Levi-Civita connexion onto $\mathrm{ker}(\theta)\cap\mathrm{ker}(\ud t)$ defines a hermitian connexion $\nabla^\star$ for $N_\star$.

That
\begin{equation*}
\nabla^\star:=\nabla-\frac{(\theta\circ\nabla)\otimes R}{\mathrm{g}(R,R)}+(\tau\circ\nabla)\otimes X
\end{equation*}is a hermitian connexion follows from it being a complex connexion, $\nabla^\star\mathrm{j}_\star =0$.

Using an arbitrary vector field $Y\in\mathrm{ker}(\theta)\cap\mathrm{ker}(\ud t)$, one can compute the curvature of $\nabla^\star$ directly from its definition to obtain
\begin{align*}
\sqrt{-1}\mathrm{curv}(\nabla^\star)
&=-\frac{\mathrm{g}(\mathrm{curv}(\nabla)Y,\mathrm{j}_\star Y)}{\mathrm{g}(Y,Y)} \nonumber \\
&\quad+\frac{\theta(\nabla Y)\wedge\theta(\nabla\mathrm{j}_\star Y)}{\mathrm{g}(R,R)\mathrm{g}(Y,Y)}+\frac{\tau(\nabla Y)\wedge\tau(\nabla\mathrm{j}_\star Y)}{\mathrm{g}(Y,Y)} \ .
\end{align*}However, if $Y$ satisfies $\mathrm{g}(Y,Y)=1$ (locally) it can be understood as a (local) unitary section of $N_\star$, and the associated potential $1$-form for $\nabla^\star$ is
\begin{equation*}
\sqrt{-1}\nabla^\star Y=-\frac{\star 1(Y,\nabla Y,R,X)}{\sqrt{\mathrm{g}(R,R)}}\otimes Y \ ;
\end{equation*}therefore, the curvature can be written (locally) as
\begin{equation*}
\sqrt{-1}\mathrm{curv}(\nabla^\star)=\ud\left(-\frac{\star 1(Y,\nabla Y,R,X)}{\sqrt{\mathrm{g}(R,R)}}\right)=-\ud\mathrm{g}(\nabla Y,\mathrm{j}_\star Y) \ .
\end{equation*}


\subsection{Normal complex line bundle}

\hspace{1.5em}The vectors fields $R$ and $X$, as well as their associated $1$-forms $\theta$ and $\tau$, can be used to define a trivial hermitian line bundle with hermitian connexion. For each point of $M$ (where $R$ does not vanish), the distribution of vector fields generated by $R$ and $X$ provides a plane subbundle of $TM$, and a complex structure on each fibre of this plane subbundle is defined by
\begin{equation*}
\mathrm{j}:=-\frac{\tau\otimes R}{\sqrt{\mathrm{g}(R,R)}}-\frac{\theta\otimes X}{\sqrt{\mathrm{g}(R,R)}} \ .
\end{equation*}

If $N$ denotes this complex line bundle, $X$ is by definition a globally defined section of it; thus, not only $N$ is trivial, but
\begin{equation*}
\mathrm{h}:=-\frac{c^2\ud\tau\circ(\mathds{1}\oplus\mathrm{j})}{\sqrt{\mathrm{g}(R,R)}}-\sqrt{-1}\frac{c^2\ud\tau}{\sqrt{\mathrm{g}(R,R)}}
\end{equation*}is a hermitian inner product, which can also be written as
\begin{equation*}
\mathrm{h}=\frac{\theta\otimes\theta}{\mathrm{g}(R,R)}+\tau\otimes\tau+\sqrt{-1}\frac{\theta\wedge\tau}{\sqrt{\mathrm{g}(R,R)}} \ .
\end{equation*}

Since $X$ can be understood as a unitary section, given any $\nu\in\Omega^1(M)$,
\begin{equation*}
\sqrt{-1}\nabla^N X:=-\nu\otimes X
\end{equation*}defines a hermitian connexion for $N$ with
\begin{equation*}
\sqrt{-1}\mathrm{curv}(\nabla^N)
=-\ud\nu \ .
\end{equation*}As a consequence, one is faced with many possibilities, but no guide to make a meaningful geometrical (or physical) choice.

An option is to mimic the potential $1$-form for $\nabla^\star$ with respect to $Y$, by recognising that it is minus the imaginary part of $\mathrm{h}_\star$ evaluated on $Y\oplus\nabla Y$. Thus, using minus the imaginary part of $\mathrm{h}$ evaluated on $X\oplus\nabla X$,
\begin{equation*}
\nu=\frac{\theta(\nabla X)}{\sqrt{\mathrm{g}(R,R)}}
\end{equation*}and
\begin{equation*}
\ud\nu=\frac{\ud\theta(\nabla X)}{\sqrt{\mathrm{g}(R,R)}}-\frac{\theta(\nabla R)\wedge\theta(\nabla X)}{\mathrm{g}(R,R)\sqrt{\mathrm{g}(R,R)}} \ .
\end{equation*}Which is a projection of $\nabla X$ onto $R$.

\begin{remark}
Projecting the Levi-Civita connexion orthogonally onto $\langle R\rangle\oplus\langle X\rangle$ will not, in general, provide a complex connexion on $N$; otherwise the Levi-Civita connexion would be complex with respect to $\mathrm{j}_\star\oplus\mathrm{j}$, implying that $TM$ is a holomorphic bundle. Also, the complex connexion
\begin{equation*}
\tfrac{1}{2}\nabla-\tfrac{1}{2}(\mathrm{j}_\star\oplus\mathrm{j})\nabla(\mathrm{j}_\star\oplus\mathrm{j})
\end{equation*}orthogonally projected onto $\langle R\rangle\oplus\langle X\rangle$ yields a flat connexion on $N$.
\end{remark}


\section{Symplectic structure}

\hspace{1.5em}Following the notation and definitions from the previous section, one is induced to conjecture a symplectic structure for (possibly not all of the points of) an orientable $4$-dimensional globally hyperbolic manifold.

\begin{conjecture}\label{symplecticstructure}
If $N_\star$ is not a trivial hermitian line bundle and $\nu$ is the proposed potential $1$-form for $\nabla^N$, then
\begin{equation*}
\varpi:=\sqrt{-1}\mathrm{curv}(\nabla^\star)+\sqrt{-1}\mathrm{curv}(\nabla^N)
\end{equation*}is a symplectic structure where $R$ does not vanish.
\end{conjecture}

The closedness of $\varpi$ is a consequence of it being the sum of curvatures of complex line bundles. Minkowski spacetime has $Q=\CR^3$ and de Sitter space has $Q=S^3$, both manifolds have trivial second cohomology group, rendering impossible to obtain a nontrivial complex line bundle over them. It is also true that $R=0$ for these examples; however, by considering a contractible subset of Minkowski spacetime and a isometry defined over it (representing the referential frame of an observer with constant proper acceleration along a spacelike line), one can make $R\neq 0$. Section \ref{schwarzschild} shows that the Schwarzschild solution has $R$ nowhere vanishing and $\varpi$ nondegenerated; there $Q=\CR^3-\{0\}$, sharing the same second cohomology group to $S^2$, which is (homeomorphic to) $\CZ$, and the de Rham class of $\varpi$ gives a generator.


\subsection{Symplectic foliation}\label{symplecticfoliationsection}

\hspace{1.5em}There is yet another Poisson structure to be considered apart from the one associated to $\varpi$, which is well defined and nontrivial (over the points where $R$ does not vanish) even when $N_\star$ is trivial. The leaves of the foliation given by $\theta$ are symplectic manifolds, since they are all complex curves; thus, it is not surprising that there is an underlying Poisson structure.

The $2$-form obtained from the imaginary part of $\mathrm{h}$
\begin{equation*}
\mathrm{Im}(\mathrm{h})=-\frac{c^2\ud\tau}{\sqrt{\mathrm{g}(R,R)}}\in\Omega^2(M) \ ,
\end{equation*}is such that
\begin{equation*}
\star\mathrm{Im}(\mathrm{h})=-\frac{\star(\tau\wedge\theta)}{\sqrt{\mathrm{g}(R,R)}}=-\frac{\imath_R\star\tau}{\sqrt{\mathrm{g}(R,R)}}=\mathrm{g}\circ(\mathrm{j}_\star\oplus\mathds{1}) \ ;
\end{equation*}i.e. $\star\mathrm{Im}(\mathrm{h})$ restricts to the area forms on the leaves given by $\theta$. In particular, both $R$ and $\depp{}{t}$ belong to its kernel, and due to $\theta\wedge\star\tau$ being zero,
\begin{equation*}
-\frac{\theta}{\sqrt{\mathrm{g}(R,R)}}\wedge\star\mathrm{Im}(\mathrm{h})=\star\tau
\end{equation*}is a volume form on $Q$ wherever $R$ does not vanish. Moreover,
\begin{equation*}
\mathrm{Im}(\mathrm{h})\wedge\star\mathrm{Im}(\mathrm{h})=\sqrt{\mathrm{g}(R,R)}\star 1
\end{equation*}defines a volume form on $M$ wherever $R$ does not vanish.

\begin{remark}
By restricting $N$ and $\star\mathrm{Im}(\mathrm{h})$ to a leaf given by $\theta$, one has its normal bundle and area form. This motivated the author's choice to use $\star$ in the symbols $N_\star$, $\mathrm{h}_\star$, $\mathrm{j}_\star$ and $\nabla^\star$.
\end{remark}

Although $\star\mathrm{Im}(\mathrm{h})$ cannot define a symplectic structure, together with $\theta$ they define a symplectic foliation on $Q$. Unfortunately, there is no consensus in the literature regarding the nomenclature of this geometric structure. Some authors call it cosymplectic, others might treat it as a symplectic vector bundle or a Poisson manifold \cite{GMP11}. In this note the relevant structure is the Poisson bracket $\{\cdot,\cdot\}_{\star}$ induced by the symplectic foliation on $M$, for which
the leaves are the complex curves with area form given by $\star\mathrm{Im}(\mathrm{h})$. The lorentzian structure allows one to construct from a function $f\in C^\infty(M)$ its gradient vector field $\mathrm{grad}(f)\in\mathfrak{X}(M)$, and applying $\mathrm{j}_\star$ to it yields
\begin{equation*}
\{f,\cdot\}_\star=\mathrm{j}_\star\mathrm{grad}(f) \ .
\end{equation*}From this expression it follows that the time function is a Casimir,
\begin{equation*}
\{t,\cdot\}_\star=0 \ ,
\end{equation*}since
\begin{equation*}
\mathrm{grad}(t)=\frac{1}{c\e^\phi}X
\end{equation*}belongs to the kernel of $\mathrm{j}_\star$.


\subsection{Quantum spacetime}\label{quantumgravity}

\hspace{1.5em}Once there is a symplectic structure for a spacetime, it is natural to consider its geometric quantisation \cite{Kostant}. Here the prequantum line bundle will be the hermitian line bundle $\mathscr{L}=N_\star\otimes N$ with hermitian connexion
\begin{equation*}
\nabla^\varpi=\nabla^\star\otimes\mathds{1}+\mathds{1}\otimes\nabla^N \ ,
\end{equation*}since it satisfies
\begin{equation*}
\sqrt{-1}\mathrm{curv}(\nabla^\varpi)=\varpi \ .
\end{equation*}Polarisations will not be considered, because the symplectic structure is defined on $M$ (as oppose to its cotangent bundle): the Hilbert space $\mathscr{H}$ will be formed by sections of the prequantum line bundle which are square integrable (with respect to the volume induced by $\varpi$), $L^2(\Gamma(\mathscr{L}))$.

To each function $f\in C^\infty(M)$, one associates its hamiltonian vector field $H_f$, the unique solution of
\begin{equation*}
\imath_{H_f}\varpi=-\ud f \ ,
\end{equation*}and the linear operator
\begin{equation*}
\hat{f}=-\sqrt{-1}\nabla^\varpi_{H_f}+f
\end{equation*}acting on $\mathscr{H}$, and if $h\in C^\infty(M)$ is another function
\begin{equation*}
\widehat{\{f,h\}}=\sqrt{-1}[\hat{f},\hat{h}] \ .
\end{equation*}

During the decade of 1970, several authors observed how to describe the interaction between the electromagnetic field and a charged particle (what in the Physics literature is known as minimal coupling) via a modification of the canonical symplectic structure of the cotangent bundle of the pertinent spacetime. This is summarised (and extended to classical Yang--Mills fields) in \cite{Sternberg}. Mimicking this idea, if $\lambda\in\Omega^1(T^*M)$ is the tautological $1$-form defined via the projection $\mathrm{pr}:T^*M\rrightarrow M$, then the relevant symplectic structure on $T^*M$ is $\ud\lambda+\kappa_G\hbar\mathrm{pr}^*(\varpi)$, where $\kappa_G$ is to be interpreted as some coupling constant. This constant (which could be a function as long as it keeps the $2$-form closed and integral) might be universal or dependent on the particle that is coupling with the gravitational field.

For the electromagnetic interaction, the $2$-form that couples with the canonical symplectic structure $\ud\lambda$ is exact, and is the tensor field that is quantised to obtain a quantum field theory. The $2$-form $\varpi$ cannot be exact (as $N_\star$ must be nontrivial), so it is possible that the correct $2$-form that couples with $\ud\lambda$ is $\sqrt{-1}\mathrm{curv}(\nabla^N)$ (which would free $\kappa_G$ from the restriction of keeping $\varpi$ integral). Moreover, this analogy suggests that the tensor field to be quantised to obtain a quantum field theory for gravity is either $\varpi$ or $\sqrt{-1}\mathrm{curv}(\nabla^N)$.

\begin{remark}
There is an embarrassing freedom regarding the actual choice. Not only one can attempt to use $\varpi$ or $\sqrt{-1}\mathrm{curv}(\nabla^N)$, but $\nabla^N$ could be defined by any $1$-form, as long as $\widehat{\{f,h\}}=\sqrt{-1}[\hat{f},\hat{h}]$ holds; even if it does not make $\varpi$ symplectic, which is the case when $\nabla^N$ is flat (but square integrable sections would need to be defined using the volume induced by $\mathrm{g}$).
\end{remark}

One would be tempted to include the mass of the particle in the coupling $\kappa_G$, as an analog for the charge of a particle in electromagnetism; however, photons do interact with gravity. An alternative would be to use the energy of the particle as an analog for the charge, or assume that the coupling is a universal constant, meaning that all particles interact in the same way with gravity. Either way, this can be tested.

Supposing $x\in C^\infty(M)$ represents a position function (e.g. a coordinate function), its quantum operator $\hat{x}$ is no longer only a multiplication, as it is in standard treatments of Quantum Mechanics. The extra term $-\sqrt{-1}\kappa_G\hbar\nabla^\varpi_{H_x}$ is a testable modification that might also introduce testable noncommutativity between different position functions (the momentum operators are not modified).

\begin{remark}
The suggested minimal coupling is not to be taken classically. The dynamics of massive particles should still be described by timelike geodesics, and light by null geodesics ---there is no modification introduced on the canonical symplectic structure. It is only when performing a quantisation of a classical system that this modified canonical symplectic structure should be considered, e.g.: if one applies geometric quantisation to a free particle moving around spacetime, one would take not $\ud\lambda$ as the curvature of the prequantum line bundle, but $\ud\lambda+\kappa_G\hbar\mathrm{pr}^*(\varpi)$, which in turn would modify the position operators (and not the momentum operators).
\end{remark}


Another approach to quantisation is to use a deformation of the Poisson algebra on $C^\infty(M)$ induced not by $\varpi$, but by the Poisson bracket $\{\cdot,\cdot\}_{\star}$ defined in subsection \ref{symplecticfoliationsection}. Although not the same, the overall result would be similar: momentum functions are kept unchanged, and different position functions would no longer commute.


\section{Schwarzschild spacetime}\label{schwarzschild}

\hspace{1.5em}Considering $Q=\CR^3-\{0\}$, on $M=Q\times\CR$ it is possible to define a lorentzian structure $\mathrm{g}$. Using linear coordinates $x,y,z,t\in C^\infty(\CR^4)$, their corresponding spherical coordinates satisfy
\begin{equation*}
\left\{\begin{array}{l}
x=r\sin(u)\cos(v) \\
y=r\sin(u)\sin(v) \\
z=r\cos(u) \\
\end{array}\right . \ ,
\end{equation*}and, for a given real number $m\in\CR$,
\begin{align*}
\mathrm{g}&=r^2\ud u\otimes\ud u+r^2\sin(u)^2\ud v\otimes\ud v  \nonumber \\
&\quad+\left(1-\frac{2Gm}{c^2r}\right)^{-1}\ud r\otimes\ud r-\left(1-\frac{2Gm}{c^2r}\right)c^2\ud t\otimes\ud t
\end{align*}is a particular solution of Einstein's equation obtained by Schwarzschild in 1916 \cite{Wald}.

By defining the function
\begin{equation*}
\phi=\ln\left(1-\frac{2Gm}{c^2r}\right)^{\tfrac{1}{2}}
\end{equation*}one can rewrite the lorentzian structure as
\begin{equation*}
\mathrm{g}_0-c^2\e^{2\phi}\ud t\otimes\ud t \ ,
\end{equation*}where $\mathrm{g}_0$ is a riemannian structure on $Q$. This is in accordance to theorem \ref{hyperbolicspacetime}, since the Schwarzschild solution is an example of a $4$-dimensional globally hyperbolic manifold \cite{Wald}.

The spacelike vector field
\begin{equation*}
R=-\frac{Gm}{r^2}\dep{}{r}
\end{equation*}can be directly obtained as the gradient of $-c^2\phi$, whilst the relevant unitary timelike vector field is
\begin{equation*}
X=-\left(1-\frac{2Gm}{c^2r}\right)^{-\tfrac{1}{2}}\frac{1}{c}\dep{}{t} \ .
\end{equation*}From the volume form
\begin{equation*}
\star 1=cr^2\sin(u)\ud u\wedge\ud v\wedge\ud r\wedge\ud t \ ,
\end{equation*}one has the $2$-form
\begin{equation*}
\omega=-\frac{1}{4\pi G}\imath_R\imath_X\star 1=\frac{m}{4\pi}\left(1-\frac{2Gm}{c^2r}\right)^{-\tfrac{1}{2}}\sin(u)\ud u\wedge\ud v \ ,
\end{equation*}which is both degenerate and not closed; for $\omega\wedge\omega=0$ and
\begin{equation*}
\ud\omega=-\frac{\mathrm{g}(R,R)}{4\pi Gc^2\e^\phi}r^2\sin(u)\ud u\wedge\ud v\wedge\ud r=-\ud\phi\wedge\omega
\end{equation*}provides a volume form on $Q$.

Because $-\ud\phi$ is an exact $1$-form (in particular it is closed), its kernel when restricted to the $3$-dimensional manifold $Q$ defines an integrable distribution (the orthogonal complement of $R$) for which its leaves are symplectic with respect to the $2$-form $\omega$. The leaves are $2$-dimensional (in this case they are the concentric spheres with centre at the origin); therefore, $\omega$ restricted to them must be closed. The nondegeneracy of the restriction of $\omega$ is justified by $-\ud\phi\wedge\omega$ defining a volume form on $Q$.

\begin{remark}
Schwarzschild's solution has that $\e^\phi\omega$ is closed, which allows for $\e^\phi(\omega+\star\omega)$ to be chosen as a symplectic structure. One might be tempted (as the author in a first try) to generalise this to obtain a symplectic structure for orientable $4$-dimensional globally hyperbolic manifolds; however, this is a particularity of Schwarzschild's solution, and in general one cannot find a function $\varphi$ for which $\e^\varphi\omega$ is closed. For example, when the lorentzian structure represents an electric charged black hole, $\e^\phi\omega$ is not closed, since its integral over any $2$-dimensional compact, without boundary, and oriented submanifold $\Sigma\subset Q$ such that $[\Sigma]\in H_2(\CR^3)$ is trivial yields a physical quantity related to the mass inside $\Sigma$; as oppose to newtonian gravity, this value depends on the representative $\Sigma$, because in General Relativity energy is also a source of gravity. The $2$-form $\mathrm{Im}(\mathrm{h})$ from subsection \ref{symplecticfoliationsection} is exactly
\begin{equation*}
-\frac{4\pi G}{\sqrt{\mathrm{g}(R,R)}}\star\omega
\end{equation*}(this holds in general), and integrating $\e^\phi\omega$ over $\Sigma$ for the Schwarzschild's solution returns the number $m\in\CR$.
\end{remark}

Now considering the Levi-Civita connexion,
\begin{equation*}
\nabla\depp{}{u}=\displaystyle-\left(1-\frac{2Gm}{c^2r}\right)r\ud u\otimes\dep{}{r}+\frac{\cos(u)}{\sin(u)}\ud v\otimes\dep{}{v}+\frac{1}{r}\ud r \otimes\dep{}{u} \ ,
\end{equation*}with respect to the unitary section $\tfrac{1}{r}\depp{}{u}$ (understood as a section of $N_\star$),
\begin{equation*}
\sqrt{-1}\nabla^\star\tfrac{1}{r}\depp{}{u}=-\cos(u)\ud v\otimes\displaystyle\frac{1}{r}\dep{}{u}
\end{equation*}and its curvature is
\begin{equation*}
\sqrt{-1}\mathrm{curv}(\nabla^\star)=\sin(u)\ud u\wedge\ud v \ .
\end{equation*}As for the hermitian line bundle $N$,
\begin{equation*}
\nabla\depp{}{t}=\displaystyle\frac{Gm}{r^2}\left(1-\frac{2Gm}{c^2r}\right)\ud t\otimes\dep{}{r}+\frac{Gm}{cr^2}\left(1-\frac{2Gm}{c^2r}\right)^{-1}\ud r\otimes\dep{}{t} \ ,
\end{equation*}
\begin{equation*}
\sqrt{-1}\nabla^N X=-\frac{Gm}{cr^2}\ud t\otimes X \ ,
\end{equation*}and
\begin{equation*}
\sqrt{-1}\mathrm{curv}(\nabla^N)=\frac{2Gm}{cr^3}\ud r\wedge\ud t \ .
\end{equation*}Therefore,
\begin{equation*}
\varpi=\sin(u)\ud u\wedge\ud v+\frac{2Gm}{cr^3}\ud r\wedge\ud t
\end{equation*}and
\begin{align*}
\varpi\wedge\varpi&=-(\mathrm{curv}(\nabla^\star)+\mathrm{curv}(\nabla^N))\wedge(\mathrm{curv}(\nabla^\star)+\mathrm{curv}(\nabla^N)) \nonumber \\
&=-2\mathrm{curv}(\nabla^\star)\wedge\mathrm{curv}(\nabla^N)=\frac{4Gm}{c^2r^5}\star 1 \ .
\end{align*}

\begin{remark}
When interpreting the classical dynamics of particles moving around a Schwarzschild black hole as one moving under the influence of a newtonian central force, the effective radial potential energy has an inverse cubic term: the radial part of $\varpi$ also has it.
\end{remark}

Computing the hamiltonian vector fields for the (spherical) coordinate functions,
\begin{equation*}
H_u=\frac{1}{\sin(u)}\dep{}{v} \ ,
\end{equation*}
\begin{equation*}
H_v=-\frac{1}{\sin(u)}\dep{}{u} \ ,
\end{equation*}
\begin{equation*}
H_r=\frac{cr^3}{2Gm}\dep{}{t} \ ,
\end{equation*}and
\begin{equation*}
H_t=-\frac{cr^3}{2Gm}\dep{}{r} \ .
\end{equation*}Hence, their Poisson brackets are:
\begin{equation*}
\{u,v\}=\varpi(H_u,H_v)=\frac{1}{\sin(u)} \ ,
\end{equation*}
\begin{equation*}
\{u,r\}=\{u,t\}=\{v,r\}=\{v,t\}=0 \ ,
\end{equation*}and
\begin{equation*}
\{r,t\}=\varpi(H_r,H_t)=\varpi(H_r,H_t)=\frac{cr^3}{2Gm} \ .
\end{equation*}

The section
\begin{equation*}
s=\frac{1}{r}\dep{}{u}\otimes X
\end{equation*}is a (local) unitary section of $\mathscr{L}=N_\star\otimes N$, and the associated potential $1$-form of $\nabla^\varpi$ is
\begin{equation*}
-\cos(u)\ud v-\frac{Gm}{cr^2}\ud t \ .
\end{equation*}By multiplying it by a function $f\in C^\infty_{\CC}(M)$ such that
\begin{equation*}
\int_M\escalar{fs}{fs}\varpi\wedge\star\varpi=\int_M\overline{f}f\varpi\wedge\star\varpi=1 \ ,
\end{equation*}one gets a (local) unitary element $\psi=fs\in\mathscr{H}$.




\begin{thebibliography}{99}


\bibitem{BeralSanchez}Antonio Bernal and Miguel S\'{a}nchez; \\ \textit{Smoothness of time functions and the metric splitting of globally hyperbolic spacetimes}; \\ Communications in Mathematical Physics, volume 257, pages 43--50 (2005). \\ \href{https://doi.org/10.1007/s00220-005-1346-1}{\tt https://doi.org/10.1007/s00220-005-1346-1}

\bibitem{Giroux}Emmanuel Giroux; \\ \textit{Convexit\'{e} en topologie de contact}; \\ Commentarii mathematici Helvetici, volume 66, issue 4, pages 637--677 (1991). \href{http://eudml.org/doc/140253}{\tt ISSN: 0010-2571; 1420-8946/e}

\bibitem{GMP11}Victor Guillemin, Eva Miranda, and Ana Rita Pires; \\ \textit{Codimension one symplectic foliations and regular Poisson structures}; \\  Bulletin of the Brazilian Mathematical Society, New Series, volume 42, pages 607–-623 (2011).
\\ \href{https://doi.org/10.1007/s00574-011-0031-6}{\tt https://doi.org/10.1007/s00574-011-0031-6}

\bibitem{Kostant}Bertram Kostant; \\ \textit{Quantization and Unitary Representations Part I: Prequantization}; \\ Lectures in modern analysis and applications III, Lecture Notes in Mathematics, volume 170 pages 87--208 (1970). \\ \href{https://doi.org/10.1007/BFb0079068}{\tt https://doi.org/10.1007/BFb0079068}

\bibitem{Sternberg}Shlomo Sternberg; \\ \textit{Minimal coupling and the symplectic mechanics of a classical particle in the presence of a Yang-Mills field}; \\ The Proceedings of the National Academy of Sciences USA, volume 74, number 12, pages 5253--5254 (1977). \\ \href{https://doi.org/10.1073/pnas.74.12.5253}{\tt https://doi.org/10.1073/pnas.74.12.5253}

\bibitem{Wald}Robert Wald; \\ \textit{General Relativity}; \\ Chicago Press (1984). ISBN 0-226-87033-2

\bibitem{Whitehead}John Whitehead; \\ \textit{The immersion of an open 3-manifold in euclidean 3-space}; \\ Proceedings of the London Mathematical Society, volumes 3--11, issue 1, pages 81--90  (1961). \href{https://doi.org/10.1112/plms/s3-11.1.81}{\tt https://doi.org/10.1112/plms/s3-11.1.81}

\end{thebibliography}
\end{document}